\documentclass[a4paper,10pt]{article}

\usepackage[english]{babel}

\usepackage{amsmath,amssymb}
\usepackage{amsthm}
\usepackage{amscd}
\usepackage{amsfonts}
\usepackage{}
\usepackage{mathrsfs}
\newcommand{\CH}{CH^{g-1}(\mathcal{J})}
\newcommand{\G}{\Gamma}
\newcommand{\C}{\mathcal{C}}
\newcommand{\Ox}{\mathcal{O}}
\newcommand{\Oc}{\mathcal{O}_C}
\newcommand{\W}{\mathcal{W}}
\newcommand{\J}{\mathcal{J}}
\newcommand{\Q}{\mathbb{Q}}
\newcommand{\X}{\mathcal{X}}
\newcommand{\Z}{\mathbb{Z}}
\newcommand{\Si}{\Sigma}
\newcommand{\co}{\mathbb{C}}
\newcommand{\Proj}{\mathbb{P}}
\newcommand{\la}{\longrightarrow}
\newcommand{\tang}{T_{S,s}}
\newcommand{\Om}{\Omega}
\newcommand{\wc}{\omega_C}
\newcommand{\oms}{\Omega_{S,s}}

\renewcommand{\bar}[1]{\overline{#1}}
    \newtheorem{lem}{Lemma}[section]
    \newtheorem{prop}[lem]{Proposition}
    \newtheorem{thm}[lem]{Theorem}

    \newtheorem{defi}[lem]{Definition}
    
    \newtheorem{rem}[lem]{Remark}

\DeclareMathOperator{\tens}{\otimes}

\begin{document}
\title{Infinitesimal invariant and vector bundles}
\date{28 may 2005}
\author{Gian Pietro Pirola and Cecilia Rizzi
\thanks{Partially supported by:
1) MIUR PRIN 2003:
  \textup{Spazi di moduli e teoria di Lie};
2) Gnsaga; 3) Far 2004 (PV):
  \textup{Variet\`{a} algebriche, calcolo algebrico, grafi orientati e topologici.} 
The second author was partially supported by a Ph.D. studentship of the Universit\`a di Pavia and by a scholarship of Politecnico di Milano.}}
\footnotetext{\emph{Mathematics Subject Classification.} Primary 14C25. Secondary 14H40, 14C15} 

\maketitle
\begin{abstract}
We study the Saito-Ikeda infinitesimal invariant of the cycle defined
by curves in their Jacobians using rank $k+1$ vector bundles. We give a criterion for which the higher cycle class map is not trivial. When $k=2$, this turns out to be strictly linked to the Petri map for vector bundles. In this case we can improve a result of Ikeda: an explicit construction on a curve of genus $g \geq 10$  shows the existence of a non trivial element in the higher Griffiths group.
\end{abstract}

\section*{Introduction}

A curve in its jacobian defines a basic algebraic cycle which is one of the most studied in algebraic geometry. In fact many geometric properties of the curve are reflected by this cycle and Abel-Jacobi mappings, normal functions and their infinitesimal invariant are the natural objects to obtain informations on the curve. 

The aim of this paper is to show that the higher infinitesimal invariant associated to the basic cycle can be interpreted using vector bundle theory. We show that the non-vanishing infinitesimal invariant is implied by the injectivity of a natural map of vector bundles, usually called Petri map, which is completely established only for line bundles on general curves but it is studied also in case of rank $n$ in the Brill-Noether theory and in the study of moduli spaces of vector bundles (see for example \cite{Mukai}, \cite{Bertram}, \cite{Teixidor}).\\   


Let $C$ be a smooth complete complex connected curve of genus
$g>2$. The choice of a base point defines the Abel-Jacobi map
$C \la J,$ where $J$ is the jacobian variety. The image of $C$ in
$J$ is an algebraic cycle $W(C)$ of codimension $g-1$. Let
$W(C)^-$ be the image of $W(C)$ under the involution $j: J \la J$,
$j(u)=-u$. Let $$Z= W(C)-W(C)^-\in CH^{g-1}(J)$$ be the
basic cycle associated to $C$.
Clearly $Z$ is homologically equivalent to zero but, for a general curve $C$, the cycle $W(C)$
is not algebraically equivalent to $W(C)^-$ (see \cite{Ceresa}). \\

This result was obtained by degeneration in \cite{Ceresa} and by infinitesimal deformations
in \cite{H} and in \cite{Co-Pi}.
Consider a family of smooth curves $f:\C \la S$ over a smooth
variety $S$ of dimension $n$ with a section of $\C$
over $S$. Let $\J$ be the family of jacobians. The algebraic cycle
$\W$ of $\J$ is defined as the image of the morphism $\C \la \J$.
Then the basic cycle $\mathcal{Z}$ given by the difference $\W -
\W^-$ is an element in $CH^{g-1}(\J)$. The higher Abel-Jacobi map
associated to $\mathcal{Z}$ defines a "normal function" (see
\cite{Griff}).
The Griffiths infinitesimal invariant was introduced as a tool to study general normal
functions in the more general context of algebraic varieties
(see also \cite{Green} and \cite{Voisin2}). The computation of the infinitesimal invariant
in our case has been carried out in \cite{Co-Pi}.\\

The structure of algebraic cycle groups on an abelian variety is
richer than in the arbitrary case. In fact we can decompose $\W =
\sum_{\nu} \W^{\nu}\in CH^{g-1}(\J) \tens \mathbb{Q}$, where
$\W^{\nu}$ lies in $CH^{g-1}_{\nu}(\J)$, that is, multiplication
by $k \in \mathbb{Z}$ acts on the cycle as $k^* \W^{\nu}
=k^{2g-2-\nu} \W^{\nu}$ (see \cite{beau}). Hence the cycle
$\mathcal{Z}$ can be written as $\sum_{i} 2 \W^{2i+1}$. It turns
out that $\W^{2i+1}$ with $i \not =0$ lies in the kernel
of the
Abel-Jacobi map (\cite{CVG}, \cite{Ikeda}).\\

In order to study the kernel of the Abel-Jacobi map,
Saito in \cite{Sa2} introduces higher normal functions and
proceeds with the definition of higher infinitesimal invariants.
For this purpose, we recall two filtrations on the Chow groups as
in \cite{Sa1}. Let $X$ be a smooth projective variety. The first
filtration $F^{\nu}CH^k(X)$ corresponds to one of those
conjectured by Bloch and Beilinson in the theory of mixed motives \cite{BB}. It has the following basic properties: 
$$ F^1CH^k(X)= CH^k(X)_{hom} \qquad  \mbox{and}  \qquad F^2CH^k(X) \subset Ker(AJ_X^k),$$ where $AJ^k_X$ denotes the higher Abel-Jacobi map. The second filtration, $Z_lF^{\nu}CH^{k}(X)$, is an ascen\-ding
filtration of $F^{\nu}CH^k(X)$ which generalizes the group of
algebraically trivial cycles on $F^{\nu}CH^{k}(X)$.

In our case, the algebraic cycle $\W^{\nu}$ lies in $F^{\nu}CH^{g-1}(\J)$. The higher
infinitesimal invariant makes it possible to determine when $\W_s^{\nu}$ is trivial in the
higher Griffiths groups,
$$Griff^{g-1,\nu}(J_s)=F^{\nu}CH^{g-1}(J_s) / F^{\nu+1}CH^{g-1}(J_s)+Z_0F^{\nu}CH^{g-1}(J_s),$$
where $J_s$ is the fiber over a general point $s\in S$. Note that, for $\nu=1$, we have the usual Griffiths group, that is, the quotient of cycles homologous to zero modulo cycles algebraically equivalent to zero.\\

The theory of the higher infinitesimal invariant in the case of
curves in their Jacobians has been developed in \cite{Ikeda}. In
particular, for $\W^{\nu} \in F^{\nu}CH^{g-1}(\J)$, Ikeda rewrites
the higher infinitesimal invariant as a linear map
$$\phi^{\nu,1}_{s}(\W^{\nu}): V^{g-1,\nu}_{s} \la \co$$ for a
general point $s\in S$, where $V^{g-1,\nu}_{s}$ is the cohomology
of the Koszul complex
$$ \bigwedge^{\nu+1} T_{S,s} \tens H^{0}(\Om_{J_s}^{\nu+2}) \la \bigwedge^{\nu} T_{S,s} \tens H^{1}(\Om^{\nu+1}_{J_s})
\la \bigwedge^{\nu-1}T_{S,s} \tens H^{2}(\Om_{J_s}^{\nu}).$$

The formula in Proposition (3.7) of \cite{Ikeda} 
generalizes \cite{Co-Pi} using some particular maps called "adjunction maps". Through the
computation of the infinitesimal invariant, in the case of special smooth plane curves,
Ikeda provides an example of a non torsion cycle $\W^{\nu}_s$ in the higher Griffiths group $Griff^{g-1,\nu}(J_s)$.\\

In this paper, we study systematically the Saito-Ikeda
infinitesimal invariant associated to cycles $\W^{\nu}$ using
vector bundles techniques. 
Here we give a brief sketch of our argument. \\

Consider the family of curves $f:\C \la S$. Denote by $C=f^{-1}(s)$ the fiber over the generic point $s\in S$.
The differential defines the following exact sequence
\begin{equation}\label{omega intro}0 \la f^*\oms^1 \la \Om^1_{\C}|_C \la \wc \la 0.\end{equation}
We can identify the conormal bundle $f^*\oms^1$ with the
rank $n$ trivial vector bundle $T^*_{S,s} \tens \Oc$.
Note that $\Om^1_{\C}|_{C}$ is a rank $(n+1)$ vector bundle with canonical determinant. Let $$\bigwedge^{n+1} H^0(\Om^1_{\C}|_{C}) \la H^0(\wc)$$ be the map given by the determinant.

Suppose that the Kodaira-Spencer map $T_{S,s}\tens H^0(\Oc)\la H^1(T_C)$ is injective.
Moreover, assume that there is a space $U$ of dimension $n+1$ in the kernel of the coboundary map $ \delta: H^0(\wc) \la T^*_{S,s} \tens H^1(\Oc)$ and set $W=\rho^{-1}(U)$ where $\rho:H^0(\Om^1_{\C}|_{C}) \la H^0(\wc)$. Consider the restriction of the determinant map
\begin{equation}\label{delta intro} \Delta_{S,s}:\bigwedge^{n+1} W \la H^0(\wc).\end{equation}
We need to consider only the contribution of the determinant map given by the total family $f$. 
For any proper subvariety $S' \subset S$ containing $s$ of dimension $m <n$, we can restrict the family $f$ to $S'$. Let $\Delta_{S',s}$ be the determinant map analogous to (\ref{delta intro}). Set $\mathcal{V}$ be the subspace of $H^0(\wc)$ generated by all the images of $\Delta_{S',s}$ when $S'$ runs over
all the proper subvarieties of $S$ with $m<n$.
Hence we define $$\beta: \bigwedge^{n+1} W \la H^0(\wc)/\mathcal{V},$$ as the composition of $\Delta_{S,s}$ with the quotient map $H^0(\wc) \la H^0(\wc)/\mathcal{V}$.
Then one can show that $\beta$ factors through $\bigwedge^{n+1}U$ to give a map
\begin{equation}\label{adj intro}\alpha: \bigwedge^{n+1} U \la H^0(\wc)/ \mathcal{V}.\end{equation}
We call this map "adjunction map".\\

To study $\alpha$, we use the language of vector bundles. Every $n$-dimensional subspace $V$ of $H^1(T_C)$ can be identified with $T_{S,s}$ for some smooth variety $S$ of dimension $n$. The natural extension associated to $V^*$ 
defines a rank $n+1$ vector bundle which plays the role of the cotangent sheaf $\Om^1_{\C}|_{C}$.
Consider the evaluation map $W \la \Om^1_{\C}|_{C}$: its kernel is given by a rank $n$ vector bundle on $C$.
The dual of this vector bundle $M$ contains all the information to compute the map $\alpha$. Then we reduce the problem to study $M$ and its properties. Since $c_1(M)=c_1(\wc)$, it is defined the determinant map
\begin{equation}\label{det intro}\bigwedge^n H^0(M) \la H^0(\bigwedge^n M) \la H^0(\wc), \end{equation} 
which is strictly linked to the map $\alpha$ as we will show in Section \ref{criterio}.  
We have the following
\newtheorem{theorem}{Theorem}
\begin{theorem}\label{main thm intro}
Let $C$ be a smooth complete curve of genus $g$. Let $M$ be a rank $n$ vector bundle on $C$ such that $c_1(M)=c_1(\wc)$ and $h^0(M^*)=0$. Assume there is a subspace $\Pi$ of $H^0(M)$ of dimension $2n+1$ generating $M$. If the map $(\ref{det intro})$ restricted to $\bigwedge^n \Pi$ 
\begin{equation}\label{hyp intro}\phi:\bigwedge^n \Pi \la H^0(\wc), \end{equation} 
is injective, then there exists a family $f:\C \la S$ of smooth curves over a smooth variety $S$ with $\dim(S)=n$ such that the adjunction map $(\ref{adj intro})$ is not trivial.
\end{theorem}

Moreover, in Section \ref{criterio}, comparing the construction of $\Om^1_{\C}|_{C}$ and $M$, we will show that the space $\Pi\subset H^0(M)$ corresponds essentially to the dual of the space $W \subset H^0(\Om^1_{\C}|_C)$ introduced before.\\

We focus our attention on the case $n=2$. Since $M$ has canonical determinant, then $M \simeq M^* \tens \wc$.  
Consider the following commutative diagram
\begin{equation*}\label{petri2}\begin{array}{ccc}
\bigwedge^2 H^0(M) & \stackrel{\phi}{\la} & H^0(\wc)\\
\rho\downarrow &  & \downarrow \tau\\
H^0(M)\tens H^0(M^* \tens \wc) & \stackrel{\mu}{\la} & H^0(M \tens M^* \tens \wc)
\end{array}\end{equation*}
where $\tau$ is the dual of the trace map $H^1(End M) \la H^1(\Oc)$, $\rho$ is the dual of the linear map $H^0(M)^*\tens H^0(M)^* \la Hom(\bigwedge^2 H^0(M), \co)$ induced by the cup product (see \cite{Mukai} p. 155) and $\mu$ is the Petri map for rank $2$ vector bundles (see \cite{Bertram}, \cite{Mukai}). \\

The injectivity of the map $\phi$ is studied in the lower genus case for stable vector bundles in \cite{Mukai}.
The main difficulty is to deal with non decomposable tensors. We study a simplified version of the problem using $M$ as a direct sum of line bundles of the same degree, that is, the semistable case. \\

We start considering a rank $2$ vector bundle $M$ on $C$ of the form $M=A \oplus (\wc \tens A^*)$, where $h^0(A)=4$ and $\deg(A)=g-1$. Such a line bundle exists for a generic smooth curve of genus $g \geq 16$. For a suitable choice of a subspace $\Pi\subset H^0(M)$, the map $\bigwedge^2 \Pi \la  H^0(A) \tens H^0(\wc \tens A^*)$ is injective; hence the injectivity of the Petri map for line bundles implies that the map (\ref{hyp intro}) satisfies the hypothesis of Theorem \ref{main thm intro}.
Moreover, a result in \cite{Colombo} assures that there exists a smooth curve $C$ of genus $g\geq 10$, which admits a theta-characteristic $A$ with $h^0(A)=4$ and such that the linear system $|A|$ defines an embedding of $C$ in $\Proj^3$. Hence we use a theta characteristic $A$ to construct $M$ as before. It is possible to find a space $\Pi$ such that the map $\bigwedge^2 \Pi \la Sym^2 H^0(A)$ is injective and to show that our curve does not lie on quadrics in $\Proj^3$. Theorem \ref{main thm intro} allows us to conclude

\begin{theorem} \label{main prop intro}
There is a family of smooth curves $\C$ of genus $g \geq 10$ on a smooth variety $S$ of dimension $2$ for which the adjunction map $\alpha$
is not trivial.
\end{theorem}

The last part of this paper is dedicated to applications. Let $\W^2$ be the cycle in $CH^{g-1}_2(\J)$ as in the decomposition of the cycle $\W$. Ikeda's theorem assures that if the adjunction map is not trivial, then also the higher infinitesimal invariant associated to this family of curves is not zero.
So there exists a non trivial element in the higher Griffiths groups $Griff^{g-1,2}(J_s)$. \\

Hence, we have the following
\begin{theorem}\label{main thm intro 2}
For general smooth curves of genus $g \geq 16$, the cycle $\W^2_s$ is a non trivial element in the higher Griffiths group $Griff^{g-1,2}(J_s)$. Moreover, there exist curves of genus $g \geq 10$ for which the correspondent cycle $\W^2_s$ is non trivial in $Griff^{g-1,2}(J_s)$.
\end{theorem}
This theorem agrees with a result of Fakhruddin: in \cite{F}, it is proved that there are curves of genus $g\geq 11$ for which the cycle $\W_s^2$ is not algebraically equivalent to zero. \\

The paper is organized as follows.\\

In Section \ref{Mdet}, we define the "$\alpha$-determinant" map
for rank $k$ vector bundles, which is a generali\-zation of the
adjunction map referred to above. Section \ref{criterio} is dedicated to prove our
main Theorem \ref{main thm intro}. In Section \ref{adj}, we define
the adjunction map $\alpha^{n+1}_{S,s}$ giving a geometric
interpretation in terms of deformation theory of curves. Moreover,
we analyze cases $n=1$ and $n=2$ using local coordinates to make a
direct computation of the adjunction map. In Section \ref{M}, we
will exhibit a rank $2$ vector bundle on smooth curves of genus
$g\geq 10$ with the properties required in Theorem
\ref{main thm intro}.
In Section \ref{inf inv}, we review some basic results about algebraic cycles and filtrations on Chow groups. We recall the definition of the higher Griffiths' infinitesimal invariant of algebraic cycles and we explain the link between $\phi^{\nu,1}_s$ and the adjunction map; at the end, we prove Theorem \ref{main thm intro 2}.\\

\textbf{Acknowledgements:} We would like to thank Alberto Collino, Atsushi Ikeda, Shuji Saito and Enrico Schlesinger for valuable suggestions and helpful discussions. We thank the Jami institute of the Johns Hopkins University of Baltimore for its support.\\

\textbf{Notation:}
In this paper, a curve $C$ is assumed to be a smooth connected complete curve defined over $\co$. All curves we consider will have genus $g>2$.

\section{The "$\alpha$-determinant" map}\label{Mdet}
Let $C$ be a smooth curve and $L$ be a line bundle on $C$.
Consider the extension given by the following short exact sequence
\begin{equation} \label{extens}0 \la L^* \la \tilde E^* \la H^1(L^*) \tens \Oc \la 0.\end{equation}
There is a bijective correspondence between equivalence classes of extensions of $H^1(L^*)\tens \Oc$ by $L^*$ and $Ext^1(H^1(L^*) \tens \Oc, L^*)$ (see for example \cite{Weibel} p. 76-77). Since we have
$$Ext^1(H^1(L^*) \tens \Oc, L^*)\simeq H^1(L^*) \tens H^1(L^*) \simeq Hom(H^1(L^*),H^1(L^*))$$ (see for example \cite{Hart} p. 243-244), then the vector bundle $\tilde E^*$ is defined by the identity map in $Hom(H^1(L^*),H^1(L^*))$. In fact, the coboundary map of (\ref{extens}) is exactly the identity.

Let $V$ be a proper subspace of $ H^1(L^*)$ of dimension $r$: we can construct the corresponding extension
$$ 0 \la L^* \la E^*_{V} \la V \tens \Oc \la 0, $$ defined through the diagram
\begin{equation} \label{aggV}
\begin{array}{ccccccccc}
0 & \la & L^* & \la & E^*_{V} & \la & V \tens \Oc & \la & 0 \\
  & & \downarrow & & \downarrow  & & \downarrow \\
0 & \la & L^* & \la & \tilde E^*& \la & H^1(L^*) \tens \Oc & \la & 0. \\
\end{array}
\end{equation}
Dualizing (\ref{aggV}), we obtain
\begin{equation} \label{aggV2}
\begin{array}{ccccccccc}
0 & \la & H^1(L^*)^* \tens \Oc & \la & \tilde E & \la & L & \la & 0 \\
  & & \downarrow & & \downarrow  & & \downarrow \\
0 & \la & V^* \tens \Oc & \la & E_{V}& \la & L & \la & 0 .\\
\end{array}
\end{equation}
The vector bundle, $E_V$, on $C$ has rank $r+1$. The composition map $\Delta^{r+1}_{V}= g \circ f$ defines a map
given by the wedge product
\begin{equation}\label{deltaV}\Delta^{r+1}: \bigwedge^{r+1} H^0(E_{V}) \stackrel{f}{\la} H^0(\bigwedge^{r+1}E_V) \stackrel{g}{\la} H^0(L). \end{equation}

Consider the long exact sequence
$$0 \la V^* \tens H^0(\Oc) \la H^0(E_V) \stackrel{\rho}{\la} H^0(L) \stackrel{\delta}{\la} V^* \tens H^1(\Oc) \la \cdots .$$
Assume that $ \dim Ker (\delta)\geq r+1$. Fix a $(r+1)$-dimensional subspace $U$ in $Ker(\delta)$. Set $\rho^{-1}(U)= W\subset H^0(E_V)$: note that $\dim(W)=2r+1$. So we can consider the restriction of the map (\ref{deltaV}) to $\bigwedge^{r+1} W$
\begin{equation}\label{deltaV2}\Delta^{r+1}_{W}: \bigwedge^{r+1} W \hookrightarrow \bigwedge^{r+1}H^0(E_V) \la H^0(L). \end{equation}

For a proper vector subspace $Z \subset V$ of dimension $s$ of $H^1(L^*)$, we construct, with the same techniques, a rank $s+1$ vector bundle $E_Z$. Note that $U$ lies in the kernel of $\delta_Z: H^0(L) \la Z^* \tens H^1(\Oc)$: then we have a subspace $W_Z$ of $H^0(E_Z)$ of dimension $s+r+1$ given by $ \rho_Z^{-1}(U)$ where $\rho_Z: H^0(E_Z) \la H^0(L)$ is defined in the above cohomology sequence. \\

The inclusion $Z \hookrightarrow  V$ induces a map of vector space  $i_{W_Z}:W \la W_Z$.
Let $i^{s+1}_{W_Z}: \bigwedge^{s+1}W \la \bigwedge^{s+1}W_Z$ be the map induced by $i_{W_Z}$ and $\Delta^{s+1}_{W_Z}$ be the map given by the wedge product $\bigwedge^{s+1} W_Z \la H^0(L)$.
We may define $\Upsilon_{Z}= \Delta^{s+1}_{W_Z} \circ i^{s+1}_{W_Z}$
$$\Upsilon_{Z}: \bigwedge^{s+1} W \stackrel{i^{s+1}_{W_Z}}{\la} \bigwedge^{s+1}W_{Z} \stackrel{\Delta^{s+1}_{W_Z}}{\la} H^0(L).$$
Let $\mathcal{V}$ be the space generated by all the images of $\Upsilon_{Z}$, where $Z$ ranges over all proper subspaces of $V$. Set $\mathcal{Q}_{L,\mathcal{V}}= H^0(L) /\mathcal{V}$.
We consider the map
$$\beta^{r+1}: \bigwedge^{r+1} W \la \mathcal{Q}_{L,\mathcal{V}}.$$
One can show that $\beta^{r+1}$ vanishes on the kernel of $\bigwedge^{r+1} \rho: \bigwedge^{r+1} W \la \bigwedge^{r+1}U$. 

\begin{defi} We call "$\alpha$-determinant" map  the following map
$$\alpha^{r+1}_{U}: \bigwedge^{r+1} U \la \mathcal{Q}_{L,\mathcal{V}}.$$ 
\end{defi}

        \section{Vector bundles and "$\alpha$-determinant"}\label{criterio}

Using vector bundles, we will give a non vanishing criterion for the "$\alpha$-determinant" map introduced in the previous section. We begin by the construction of the rank $k+1$ vector bundle $E_V$ starting from a fixed rank $k$ vector bundle $M$, with same special properties. Later on, we will show how conditions imposed on $M$ and on its sections imply the non triviality of "$\alpha$-determinant" map. \\

Let $L$ be a line bundle on $C$ with $\deg(L)>0$.
From now on, let $M$ be a rank $k$ vector bundle on $C$ with $c_1(M)=c_1(L)$ and $h^0(M^*)=0$. Assume that $\Pi$ is a $(2k+1)$-dimensional subspace of $H^0(M)$ and $\Si$ is a $(k+1)$-dimensional subspace of $\Pi$, which generates $M$. 
Look at the evaluation map $\Si \tens \Oc \la M$; we have a short exact sequence  
\begin{equation}\label{succ M*} 0 \la N \la \Si \tens \Oc \la M \la 0,\end{equation}
where $N$ is a line bundle on $C$. Computing the first Chern class, we obtain that $c_1(N)+c_1(M)=0$.   
Since $c_1(M)= c_1(L)$, it follows that $$N=L^*.$$ 
Consider the long exact sequence in cohomology
\begin{equation*}
0 \la \Si \tens H^0(\Oc) \stackrel{\psi }{\la} H^0(M) \stackrel{\delta}{\la} H^1(L^*) \stackrel{\varphi}{\la} \Si \tens H^1(\Oc) \la H^1(M) \la 0.
\end{equation*}
\noindent
Then $\Si$ lies in the kernel of $\delta$ and hence, there is a subspace $V=\Pi/\Si$ of $H^1(L^*)$ of dimension $k$: by abuse of notation, we also denote by $V$ the subspace of $H^0(M)$ whose elements, via $\delta$, generate $Im (\delta)$. \\

\noindent
Let $E_V^*$ be the kernel of the evaluation map $\Pi \tens \Oc \la M$.
We construct the following diagram
\begin{equation}\label{diag1}
\begin{array}{ccccccccc} & & 0 & & 0 & & \\
& & \uparrow  & & \uparrow  & & \\
0 & \la & M & \la & M & \la & 0 &  &  \\
& & \uparrow  & & \uparrow  & & \uparrow \\
0 & \la & \Si \tens \Oc & \la & \Pi\tens \Oc & \la & V\tens \Oc & \la & 0 \\
& & \uparrow  & & \uparrow  & & \uparrow \\
0 & \la & L^* & \la & E_V^* & \la & V \tens \Oc& \la & 0\\
& & \uparrow  & & \uparrow & & \uparrow \\
& & 0 & & 0 & & 0.
\end{array}
\end{equation}

\noindent
Hence we obtain the exact sequence
$$0 \la L^* \la E_V^* \la V \tens \Oc \la 0.$$
Therefore, we summarize the construction in the following
\begin{prop}\label{teo1}
Let $M$ be a rank $k$ vector bundle on $C$ such that $c_1(M)=c_1(L)$ and $h^0(M^*)=0$. Assume that $\Pi$ is a subspace of $H^0(M)$ of dimension $2k+1$ which generates $M$ and $\Si$ is a subspace of $\Pi$ with $dim(\Si)=k+1$ such that it generates $M$. Set $V=\Pi/\Si$. Then there exists a rank $k+1$ vector bundle $E_V$ on $C$ with at least $2k+1$ non trivial sections and there is the following exact sequence
\begin{equation}\label{E^*_V} 0 \la L^* \la E_V^* \la V \tens \Oc \la 0.\end{equation}
\end{prop}
The sequence (\ref{E^*_V}) is exactly that we have considered in section \ref{Mdet}.
Consider the following diagram induced in cohomology by the dual of (\ref{diag1})
\begin{equation*}\label{cohom}
\begin{array}{cccccccccc} 
0 & \la & V^* \tens H^0(\Oc) & \la & \Pi^* \tens H^0(\Oc) & \la & \Si^* \tens H^0(\Oc) & \la & V^* \tens H^1(\Oc)& \la \\
& & \downarrow  & & \downarrow  & & \downarrow & & \downarrow \\
0 & \la & V^* \tens H^0(\Oc)& \la & H^0(E_V) & \la & H^0(L) & \la & V^* \tens H^1(\Oc)& \la 
\end{array}
\end{equation*}
Hence, looking at the construction of $E_V$, we can identify the space $\Pi^*$ with $W$, where $W$ is the space of sections of $E_V$ introduced in the previous section. Then, for convenience of notations, in the rest of this section, we set $\Pi=W^*$. \\
Consider the composition map $\phi$ 
\begin{equation*}\label{phi}\phi: \bigwedge^k W^* \hookrightarrow  \bigwedge^k H^0(M) \la H^0(L).\end{equation*}
given by the restriction to $\bigwedge^k W^*$ of the determinant map $\bigwedge^k H^0(M) \la H^0(\bigwedge^k M) \la H^0(L)$.\\

The condition $h^0(M^*)=0$ implies that there is an inclusion $\Si \la H^0(L)$ in the cohomology sequence  of the dual of (\ref{succ M*}). Let $U$ be the image of this inclusion in $H^0(L)$. Consider the dual sequence of (\ref{E^*_V}) $$0 \la V^* \tens \Oc \la E_V \la L \la 0.$$
Note that $U$ is a $(k+1)$-dimensional subspace of $H^0(L)$ which lies in the kernel of the coboundary map $\delta: H^0(L) \la V^* \tens H^1(\Oc)$. So we can consider $\rho^{-1}(U) \subset H^0(E_V)$ and, using the dual diagram of (\ref{diag1}), we can identify $\rho^{-1}(U)= W$.\\

We give a criterion to establish whether the "$\alpha$-determinant" map $\alpha_U^{k+1}: \bigwedge^{k+1} U \la \mathcal{Q}_{L,\mathcal{V}}$ is zero. 

\begin{thm}\label{teo2} Under hypotheses of Proposition (\ref{teo1}), suppose that $M$ is the rank $k$ vector bundle and $W^*$ is the subspace of $H^0(M)$ introduced before. If the map $$\phi:\bigwedge^k W^* \la H^0(L)$$ is injective,
then the "$\alpha$-determinant" map $\alpha^{k+1}_{U}: \bigwedge^{k+1} U \la \mathcal{Q}_{L,\mathcal{V}}$
is not zero.
\end{thm}
\begin{proof}
We begin by noting that the injectivity of the map $\Si \la U $ implies that also the map $\bigwedge^{k+1} \Si \la \bigwedge^{k+1} U $ is injective. We want to show that the map $\alpha^{k+1}_U: \bigwedge^{k+1} U \la \mathcal{Q}_{L,\mathcal{V}}$ is not zero; actually we will prove that the map $\bigwedge^{k+1} \Si \la \mathcal{Q}_{L,\mathcal{V}}$ is not zero. \\

Note that $\dim(W^*)=2k+1$; then we have an isomorphism $\bigwedge^{k} W^* \simeq \bigwedge^{k+1} W$ (defined up to a constant). Then we can rewrite $\phi$ as the injection $\phi': \bigwedge^{k+1}W \la H^0(L)$.\\

Let $\mathcal{H}$ be the kernel of the map $\bigwedge^{k+1} W \la \bigwedge^{k+1} \Si$.
We have to show that the map $\chi$ defined in the following diagram is surjective
$$
\begin{array}{ccccccccc}
0 & \la & \mathcal{H} & \la & \bigwedge^{k+1} W & \la & \bigwedge^{k+1} \Si & \la & 0\\
& & {\chi }\;{\downarrow}\;\;\;  & & {\phi'}\;{\downarrow}  & & {\alpha^{k+1}_U}\;{\downarrow}\;\;\;\;\;\; \\
0 & \la & \mathcal{V} & \la & H^0(L) & \la & \mathcal{Q}_{L,\mathcal{V}} & \la & 0.\\
\end{array} $$
The space $\mathcal{V}$ was defined in section \ref{Mdet} as the subspace of $H^0(L)$ spanned by the image of $\Upsilon_Z$,
where $Z$ is a proper subspace of $V$ of dimension $r$. Following the construction of Proposition (\ref{teo1}), let $E_Z$ be the rank $r+1$ vector bundle associated to $Z$ and $W_Z$ be a subspace of $H^0(E_Z)$ of dimension $k+r+1$ given by $\rho_Z^{-1}(U)$.\\

Let $Y^*$ be the kernel of the map $V^* \la Z^*$: note that $\dim(Y^*)=k-r$.  Then we have the exact sequence
\begin{equation}\label{Y} 0 \la Y^* \la W \la W_Z \la 0 .\end{equation}

\noindent
Consider the $(k+1)$-wedge product of the sequence before
$$0 \la W_{V,Z} \la \bigwedge^{k+1} W \la \bigwedge^{k+1} W_Z \la 0;$$
the kernel $W_{V,Z}$ is given by $\bigoplus_{p=r}^{k-1}( \bigwedge^{p+1} W_Z \tens \bigwedge^{k-p} Y^*)$.
Let $\mathcal{H}_Z$ be the kernel of $\bigwedge^{k+1} W_Z \la \bigwedge^{k+1} \Si$.
We construct a new commutative diagram
\begin{equation*}\label{diag4}
\begin{array}{ccccccccc} & & & & 0 & & 0\\
& &   & & \uparrow  & & \uparrow \\
& & 0 & \la & \bigwedge^{k+1} \Si & \la & \bigwedge^{k+1} \Si &\la & 0\\
& & \uparrow  & & \uparrow  & & \uparrow \\
0 & \la & W_{V,Z} & \la & \bigwedge^{k+1} W & \la & \bigwedge^{k+1} W_Z & \la & 0 \\
& & \uparrow  & & \uparrow  & & \uparrow \\
0 & \la & W_{V,Z} & \la & \mathcal{H} & \la & \mathcal{H}_{Z} & \la & 0 \\
& & \uparrow & & \uparrow  & & \uparrow \\
& & 0 & & 0 & & 0.
\end{array}
\end{equation*}
Consider $\bigwedge^{r+1} W_Z \tens \bigwedge^{k-r} Y^* \subset W_{V,Z}$.
Now fixing a basis $y_1,\cdots, y_{k-r}$ of the space $Y^*$,
we have the isomorphism $\bigwedge^{r+1} W_Z \simeq \bigwedge^{r+1} W_Z \tens \bigwedge^{k-r} Y^*.$  Then $\bigwedge^{r+1}W_Z$ lies in $\mathcal{H}$.\\

\noindent
It remains to show that the images in $H^0(L)$ of the following two maps
$$\begin{array}{c}\bigwedge^{r+1} W \la \bigwedge^{r+1}W_Z \la H^0(L) \\
\\
\bigwedge^{r+1} W_Z \hookrightarrow  \bigwedge^{k+1} W \la H^0(L)\end{array}$$ 
are equal. First, let $s_1,\cdots, s_{r+1}$ be $r+1$ elements in $W_Z$. Set $\tilde s_i$ be liftings of $s_i$ in $W$ for $i=1,\cdots, r+1$. Hence, the element $\tilde s_1 \wedge \cdots \wedge \tilde s_{r+1} \in \bigwedge^{r+1} W $ corresponds exactly to $s_1 \wedge \cdots \wedge s_{r+1} \in \bigwedge^{r+1} W_Z $. Locally, by a choice of the basis of $W$, the map $\bigwedge^{r+1}W \la H^0(L)$ is given by the determinant $\det[ s_1, \cdots, s_{r+1}]$.\\

\noindent
Moreover, consider the sequence (\ref{Y}): we can identify $W= W_Z \oplus Y^*$.
Through a suitable choice of a basis of $W$, that is, completing $y_1, \cdots , y_{k-r}$ with a basis of $W_Z$, the element  $s_1 \wedge \cdots \wedge s_{r+1} \wedge y_1 \wedge \cdots \wedge y_{k-r} \in \bigwedge^{r+1}W_Z \tens \bigwedge^{k-r}Y^* \subset \bigwedge^{k+1} W$ goes to 
$$\det \left[ \begin{array}{c|c} s_1, \cdots, s_{r+1} & 0 \\
\hline
\cdots & I_{k-r} \end{array}\right],$$
as a local computation of the determinant map shows. This allows us to conclude that the map $\chi$ is surjective.

\end{proof}

        \section{The adjunction map}\label{adj}

In this section, we consider the case $L=\wc$. We will give a geometric interpretation of the "$\alpha$-determinant" map using the language of infinitesimal deformations of the curve $C$. \\

Let $f:\C \la S$ be a family of smooth curves of genus $g$ over $S$, a smooth irreducible variety of dimension $n$. Consider a general point $s \in S$: we call $C$ the curve $f^{-1}(s)$ over the point $s$. If we consider the restriction at $C$ of the morphism of vector bundle $T_{\C} \la f^*(T_S)$, there is an exact sequence
\begin{equation} \label{tang}
0 \la T_C \la T_{\C}|_C \la f^*T_{S,s} \la 0.
\end{equation}
Note that $f^*T_{S,s}$ is a trivial vector bundle with fibers $T_{S,s}$. The extension of (\ref{tang}) is classified by the Kodaira-Spencer map $$\rho : T_{S,s}=H^0(C,f^*T_{S,s}) \la H^1(C,T_C),$$
induced by the long exact sequence in cohomology.
Suppose that $\rho$ is injective. By abuse of notation, we denote again with $T_{S,s}$ the image of $\rho$ in $H^1(T_C)$.
In particular, it will be more useful the dual of (\ref{tang})
\begin{equation}\label{omega} 0 \la T^*_{S,s} \tens \Oc \la \Om^1_{\C}|_{C} \la \wc \la 0 .\end{equation}
So we rewrite the map (\ref{deltaV}) given by the determinant in this way
\begin{equation}\label{adj1} \Delta^{n+1}_{S,s}:\bigwedge^{n+1} H^0(\Om^1_{\C}|_{C}) \la H^0(\bigwedge^{n+1} \Om^1_{\C}|_{C}) \la H^0(\wc).\end{equation}

The coboundary map $\delta_{[\xi]}: H^0(\wc) \la T_{S,s}^* \tens H^1(\Oc)$ is given by the cup product with the extension class of the sequence (\ref{tang}) $[ \xi ]\in T_{S,s} \tens H^1(\Oc)$. Fix $U$ be a subspace of dimension $n+1$ of the kernel of $\delta_{[\xi]}$. We can pick $W$ in $H^0(\Om^1_{\C}|_{C})$ defined by $\rho^{-1}(U)$ with $\rho:H^0(\Om^1_{\C}|_{C}) \la H^0(\wc)$. Remark that $w_i \in U $ is equivalent to ask that $\xi_j \cdot w_i=0$ with $\xi_j \in H^1(T_C)$ for $j=1,\cdots n$. \\

\noindent
Let $S'$ be a smooth irreducible proper subvariety of $S$ of dimension $m<n$, containing $s$. Consider the restriction of the family $f':\C\la S'$ over the subvariety $S'$. We can repeat the same argument for the family $f'$: in particular, we find a subspace $W'$ of dimension $m+n+1$ of $H^0(\Om^1_{\C/S'}|_{C})$.
This allows us to define, as in section \ref{Mdet}, a map
\begin{equation*} \Upsilon_{S,S',s}: \bigwedge^{m+1} W \la \bigwedge^{m+1} W' \stackrel{\Delta^{m+1}_{S',s}}{\la} H^0(\wc).\end{equation*}
where the map $\Delta^{m+1}_{S',s}$ is induced by the determinant
$\Delta^{m+1}_{S',s}: \bigwedge^{m+1} W' \la H^0(\wc).$
We denote by $\mathcal{V}$ the subspace of $H^0(\wc)$ generated by all the image of the map $\Upsilon_{S,S',s}$, where $S'$ ranges over all the proper subvarieties of $S$ with $0 \leq dim(S')=m <n$.
Notice that if $m=0$, the subvariety $S'$ is reduced to a point $\{s\}$: in this case, we are considering the space generated by the image of the map $H^0(\Om^1_{\C/\{s\}}|_C)\simeq H^0(\wc) \la H^0(\wc)$, i.e. by the image of the $n$ $1$-forms themselves. From now on, we denote $\mathcal{Q}_{\wc,\mathcal{V}}= H^0(\wc) / \mathcal{V}$.\\

\begin{rem} Let $\tilde S$ be a smooth proper subvariety of $S$ of dimension $k$. The dimension of the image of the map $\Upsilon_{S,\tilde S,s}$ in $H^0(\wc)$ is
$ \varrho_k=\frac{n+1}{k+1} \left[{n \choose k}  \right]^2.$
When $\tilde S$ ranges over all the proper subvarieties of $S$ and all these conditions are independent, as one expects, then we have $g \geq \sum_{k=0}^{n-1} \varrho_k.$
\end{rem}

Finally, we define the "$\alpha$-determinant" map which, in this case, is called \emph{adjunction map}
\begin{equation}\label{aggiun}
\alpha^{n+1}_{S,s}: \bigwedge^{n+1} U \la \mathcal{Q}_{\wc,\mathcal{V}}.
\end{equation}

Now we can prove \textbf{Theorem \ref{main thm intro}}.

\begin{proof}(\textbf{Theorem 1})\\
Let $V$ be a proper subspace of $H^1(T_C)$ of dimension $n$. Since $H^1(T_C)$ parametrizes the first order deformation, we can identify $V$ with a subspace $T_{S,s}$ of $H^1(T_C)$ for any variety $S$ of dimension $n$. Then we are able to construct the correspondent extension
\begin{equation*}\label{E^* e TC} 0 \la T_C  \la E^* \la T_{S,s} \tens \Oc \la 0.\end{equation*}
We look at the dual sequence.
\begin{equation}\label{E e TC} 0 \la T^*_{S,s} \tens \Oc \la E \la \wc \la 0.\end{equation}
The vector bundle $E$ has rank $n+1$.
The sequence (\ref{E e TC}) allows us to identify the rank $n+1$ vector bundle $E$ with the cotangent sheaf of a family $f:\C \la S$ of smooth curves of genus $g$ over $S$ restricted to a fiber $C$, $\Om^1_{\C}|_{C}$, where $S$ is the smooth irreducible variety of dimension $n$ introduced before.
Hence, to conclude the proof, it is enough to apply Proposition (\ref{teo1}) and Theorem (\ref{teo2}) in the case of $L= \wc$.
\end{proof}

In the following paragraph, we analyze in details cases $n=1$ and $n=2$. In particular we give a description of the adjunction map $\alpha^{3}_{S,s}$ using local coordinates.\\

    \subsection{Some examples: the cases $n=1$ and $n=2$}\label{n=2}

Let $\C \la S$ be a family of smooth curves of genus $g$ over a smooth variety $S$ of dimension $n$. 
We have the following exact sequence
\begin{equation}\label{n1}
0 \la N^* \la \Om^1_{\C}|_{C} \la \wc \la 0
\end{equation} where $N^*$ is the conormal bundle. Remark that if we choose a basis for the tangent space $T_{S,s}$, we have $N \simeq T_{S,s} \tens \Oc$ with the notations used before. Fix $U$ in $Ker(\delta)$ of dimension $n+1$ where $\delta$ is 
the coboundary map $ H^0(\wc) \la T^*_{S,s} \tens H^1(\Oc)$. \\

The case $n=1$ is treated by Collino and Pirola in \cite{Co-Pi}. 
The adjunction map $\alpha^2_{S,s}: \bigwedge^2 U \la \mathcal{Q}_{\wc,\mathcal{V}}$
can be described in terms of coordinates in the following way. Let $w_1,w_2 \in U$: there exist two liftings $\tilde w_1,\tilde w_2 \in H^0(\Om^1_{\C}|_{C})$. Then we have
$$\alpha^2_{S,s}(w_1 \wedge w_2)= det \left[ \begin{array} {cc} w_1 & w_2 \\ h_1 & h_2 \end{array}\right],$$
where we write $\tilde w_i= w_i dz + h_i d\varepsilon$ for $i=1,2$, with $h_i \in \C^{\infty}(C)$.
It is easy to check that the map $\alpha^2_{S,s}(w_1 \wedge w_2)$ is well defined modulo the space $\mathcal{V}$, which in this case is spanned only by $w_1$ and $w_2$.\\

Consider now a family of smooth curves over a variety $S$ of dimension $2$; the adjunction map
\begin{equation}\label{aggiun3}
\alpha^3_{S,s}: \bigwedge^3 U \la \mathcal{Q}_{\wc,\mathcal{V}},
\end{equation}
can be expressed in terms of local coordinates as follows
$$ \alpha^3_{S,s}(w_1 \wedge w_2 \wedge w_3)= det \left[ \begin{array}{ccc} w_1 & w_2 & w_3 \\ h_1 & h_2 & h_3 \\ g_1 & g_2 & g_3 \end{array} \right],$$
where $w_1,w_2,w_3 \in U$ and we write liftings $\tilde w_i \in H^0(\Om^1_{\C}|_C)$ in local coordinates $\tilde w_i= w_i dz +h_i d\varepsilon_1 + g_i d\varepsilon_2$ for $i=1,2,3$.
Observe that, in this case, $\mathcal{V}$ is the $10$-dimensional space spanned by $w_1,w_2,w_3$ and by $\alpha^2_{S,s}(w_h \wedge w_k)$ for $h,k=1,2,3$ but $h\not=k$.

    \section{Construction of a rank $2$ vector bundle}\label{M}

In this paragraph, we consider $k=2$ and $L= \wc$. We are going to exhibit an explicit rank $2$ vector bundle $M$ on a smooth curve $C$ of genus $g\geq 10$, with a space $\Pi\subset H^0(M)$, which satisfies all properties required in the hypotheses of Theorem \ref{main thm intro}. This allows us to conclude that the adjunction map $\alpha^3_{S,s}: \bigwedge^3 U \la \mathcal{Q}_{\wc,\mathcal{V}}$ is not trivial.\\

For convenience, we will indicate with $(\star )$ the whole of the following properties:
\begin{enumerate}
\item [1)] $M$ is generated by global sections;
\item [2)] $M$ has canonical determinant;
\item [3)] $h^0(M)\geq 5$;
\item [4)] $h^0(M^*)=0$;
\item [5)] there is a $5$-dimensional subspace $\Pi$ of $H^0(M)$ such that the map
$ \bigwedge^2 \Pi \la H^0(\wc)$
is injective.
\end{enumerate}

\begin{thm}\label{costruzione M}
Let $C$ be a smooth complex curve of genus $g$.
\begin{itemize}
\item[\textbf(i)] For a general smooth curve $C$ of genus $g\geq 16$, there is a rank $2$ vector bundle $M$ on $C$ which satisfies the condition $(\star)$.
\item[\textbf(ii)] For $g\geq 10$, there exists a smooth curve $C$ such that admits a rank $2$ vector bundle $M$ which satisfies the condition $(\star)$.
\end{itemize}
\end{thm}

\begin{proof} \textbf{Part (i)}

We construct $M$ as a sum of two line bundles
$$M=A \oplus (\wc \tens A^*),$$
where $A$ is a base points free line bundle with $h^0(A)=4$ and $\deg(A)=g-1$.
Observe that $A \in W^{3}_{g-1}=\{A \in Pic^{g-1}(C) : h^0(A) \geq 4\}$. Brill-Noether theory assures that the general curve of genus $g\geq 16$ admits such a linear system. It is easy to check that $M$ is a rank $2$ vector bundle on $C$ which satisfies the first four properties of ($\star$).\\

Then it remains to prove that there is a $5$-dimensional subspace $\Pi$ of $H^0(M)$ such that
$\bigwedge^2 \Pi \la H^0(\wc)$ is injective.\\

Observe that
$\bigwedge^2H^0(M)=\bigwedge^2H^0(A) \oplus \bigwedge^2 H^0(\wc \tens A^*) \oplus H^0(A) \tens H^0(\wc \tens A^*)$.
Consider the Petri map for line bundles
$$\mu: H^0(A) \tens H^0(\wc \tens A^*) \la H^0(\wc).$$
A classical result (see \cite{Gieseker}) of the theory of the moduli space of curves of genus $g$ assures that
if $[C] \in \mathcal{M}_g$ is general, the Petri map is injective for every line bundle $A$ on $C$.

Let $s_i$ and $t_i$ for $i=1,\cdots,4$ be non trivial sections of $A$ and $\wc \tens A^*$ respectively. Define $\Pi$ to be the vector space spanned by
$$ (s_1, 0), (0, t_1), (s_2, t_2), (s_3, t_3), (s_4, t_4).$$
So it remains to show that $\bigwedge^2 \Pi \stackrel{\varphi}{\la} H^0(A) \tens H^0(\wc \tens A^*)$ is injective. But this follows by a direct computation of the image of $\varphi$.
\end{proof}

\begin{proof} \textbf{Part (ii)}

As in the previous case, we want to construct $M$ as an extension of line bundles. We define $$M=A \oplus A,$$ where $A$ is a theta-characteristic such that $h^0(A)=4$. Obviously, we have $\deg(A)=g-1$. Then $M$ is a rank $2$ vector bundle on $C$ which satisfies 1-4 in $(\star)$. \\

We construct directly an example of smooth curve $C$ of genus $10$ with a theta characteristic $A$  with $h^0(A)=4$ without base points. Noting that the curve must have degree $9$, we take $C$ as complete intersection of two cubics in $\Proj^3$. The linear system $|A|$ defines an embedding of $C$ in $\Proj^3$ and so $\Oc(1)=A$. Moreover, since $A$ has to be a theta characteristic, it holds $\Oc(2)=\wc$. \\

Now look at smooth curves of genus $g \geq 11$.
The existence of such a theta-characteristic $A$ on $C$ is assured by the following
\begin{thm} [\cite {Colombo}]\label{thm colombo}
Let $\mathcal{M}_g^3$ be the subset of the moduli space $\mathcal{M}_g$ parametrizing smooth curves of genus $g$ with an even theta-characteristic whose space of sections has dimension at least $4$.
For $g\geq 9$, there is a component of $\mathcal{M}_g^3$ whose generic point corresponds to a curve $C$ with theta-characteristic $A$ such that $|A|$ defines an embedding of $C$ in $\Proj^3$.
\end{thm}

To conclude the proof, we will construct a $5$-dimensional subspace $\Pi$ of $H^0(M)$ such that $\bigwedge^2 \Pi \la H^0(\wc)$ is injective. \\

Consider the vector space $\Pi$ spanned by
$$ (s_1, 0), (0, s_2), (s_2, s_3), (s_3, s_4), (s_4, s_1)$$
where $s_i$ for $i=1,\cdots,4$ are non trivial independent sections of $A$.
\noindent
The map $\bigwedge^2 \Pi \stackrel{\varphi}{\la} Sym^2H^0(A) $ is injective: this follows by a direct computation of the image of $\varphi$

$$
\begin{array}{c}
\!\!\!\!\!\!\!\!\!\!\!\!\!\!\!\!\!\!\!\!\!\!\begin{array}{rcl}
\varphi ((s_1, 0) \wedge (0, s_2))&=& s_1 \tens s_2
\end{array}
\\
\begin{array}{rclcrcl}
\varphi ((s_1, 0) \wedge (s_2, s_3))&=&s_1 \tens s_3 &\qquad& \varphi ((0, s_2) \wedge (s_3, s_4))&=&s_2 \tens s_3 \\
\varphi ((s_1, 0) \wedge (s_4, s_1))&=&s_1 \tens s_1 &\qquad& \varphi ((0, s_2) \wedge (s_2, s_3))&=&s_2 \tens s_2 \\
\varphi ((s_1, 0) \wedge (s_3, s_4))&=&s_1 \tens s_4 &\qquad& \varphi ((0, s_2) \wedge (s_4, s_1))&=&s_2 \tens s_4 \\
\end{array}
\\
\begin{array}{rcl}
\varphi ((s_2, s_3) \wedge (s_3, s_4))&=& s_2 \tens s_4 - s_3 \tens s_3\\
\varphi ((s_2, s_3) \wedge (s_4, s_1))&=& s_2 \tens s_1 - s_3 \tens s_4\\
\varphi ((s_3, s_4) \wedge (s_4, s_1))&=& s_3 \tens s_1 - s_4 \tens s_4.
\end{array}
\end{array}
$$

Consider the embedding of the curve $C$ given by the linear system $|A|$
$$\varphi_{|A|}: C \la \Proj^3.$$ Since the kernel of the multiplication map $Sym^2H^0(A) \la H^0(\wc)$ is given by $H^0(\mathcal{I}_C(2))$, it is enough to show that $C$ does not lie on a quadric surface in $\Proj^3$.\\

It is clear that the curve $C$ complete intersection of two cubics does not lie on a quadric in $\Proj^3$.\\

Consider now curves of genus $g\geq 11$. If $C$ lies on a smooth quadric surface $Q_4 \subset \Proj^3$, we can describe it through its bidegree $(a,b)$. By the computation of the degree and the genus, we obtain that $ab-2a-2b=0$: this equality shows that such a curve of genus $g\geq 11$ does not exist on a smooth quadric surface.\\

Now it remains to show that $C$ does not lie in $Q_3$, where $Q_3$ is a rank $3$ quadric in $\Proj^3$, that is, a cone over a smooth conic $\Gamma\subset \Proj^2$, with vertex $P$. Blowing up the point $P$, we obtain a ruled surface $S=\mathbb{F}_2$. Then a curve $C$ on this surface can be written as $a C_0 +b f$ where $f$ is a fibre and $C_0$ is the section corresponding to the vertex such that $C_0^2=-2$.
The computation of the genus and the degree of $C$ establishes that $ b=a+2+\frac{4}{a-2}$. Then also this case is impossible.

\end{proof}

Here we prove \textbf{Theorem \ref{main prop intro}}.
\begin{proof}(\textbf{Theorem \ref{main prop intro}})\\
The proof is a direct consequence of Theorems (\ref{teo2}) and (\ref{costruzione M}). In fact we have seen that
on the general smooth curve $C$ of genus $g\geq 16$, there exists a rank $3$ vector bundle $E$, which can be identified with $\Om^1_{\C}|_C$ for a family $\C$ of curves of genus $g\geq 16$, such that the correspondent adjunction map $\alpha_{S,s}^3$ is not trivial. Similarly, for $g\geq 10$ there exists a smooth curve $C$ which admits a rank $3$ vector bundle $E\simeq \Om^1_{\C}|_C$ such that the map $\alpha^3_{S,s}$ is not zero where $\C$ is a family of curves of genus $g\geq 10$. \\
\end{proof}

        \section{Application to the infinitesimal invariant}\label{inf inv}

We will study the infinitesimal invariant for a family of algebraic cycles on jacobian varieties, using vector bundles. In fact, the infinitesimal invariant can be computed through the adjunction map, as Ikeda showed \cite{Ikeda}.

        \subsection{Filtrations on Chow groups}

We recall some definitions about filtrations on Chow groups, in \cite{Sa1} and \cite{Sa2}.\\

Let $X$ and $Y$ be two smooth projective varieties with $n=dim(X)$ and $m=dim(Y)$.
A \textit{correspondence} between $X$ and $Y$ is an algebraic cycle $\G \in CH^r(Y \times X)$.
It induces a map $$\G_*: CH^k(X) \la CH^{k+r-m}(Y)$$ by defining $\G_*(\alpha)= (pr_Y)_*((pr_X)^*(\alpha) \cdot \G)$,
where $pr_X: Y \times_S X \la X$ and $pr_Y: Y \times_S X \la Y$ are the projections and $\cdot$ is the product intersection of cycles.

\begin{defi}
Let $\X$ be a family of smooth projective varieties over a smooth variety $S$. We define a decreasing filtration on $CH^r(\X)$
$$F^0CH^r(\X)\supset F^1CH^r(\X) \supset \cdots \supset F^{\nu}CH^r(\X) \supset \cdots$$
in the following inductive way.
For $\nu=0$ $F^0CH^r(\X)=CH^r(\X)$;
$$F^{\nu+1}CH^r(\X)= \sum_{\mathcal{V},q,\G} Im\{ \G_{\ast}:F^{\nu}CH^{r+d_V-q}(\mathcal{V}) \la CH^r(\X)\}$$
where $\mathcal{V}$, $q$ and $\G$ range over the following data:
\begin{enumerate}
\item $\mathcal{V}$ is a family of smooth projective varieties over $S$ of dimension $d_V$;
\item $q$ is an integer such that $r \leq q\leq r+d_V$;
\item $\G \in CH^q(\mathcal{V} \times \X)$ is an algebraic cycle such that for any
$s \in S$, the map $$ \G_*: H^{2d_V-2q+2r-\nu}(V_s) \la H^{2r-\nu}(X_s)/F^{r-\nu+1}(X_s)$$ is zero.

\end{enumerate}
\end{defi}
\noindent

\begin{rem}
\textup{We can describe completely only some of these spaces $F^{\nu}CH^r(\X)$.
We know that $$F^1CH^r(\X)=CH^r(\X)_{hom},$$
where $CH^r(\X)_{hom}$ denotes the subgroup of cycle classes which are homologically equivalent to zero. Moreover, we know that $$F^2CH^r(\X) \subset Ker(AJ^r_\X),$$ where $AJ^r_\X:CH^r(\X)_{hom} \la \J^r(\X)$ is the Abel-Jacobi map.}
\end{rem}

It is possible to introduce an ascending filtration $Z_l F^{\nu}CH^{r}(\X)$ on $F^{\nu}CH^{r}(\X)$. We report only that part we need for the definition of the higher Griffiths group and refer to \cite{Sa1} for the complete subject.
We define  $Z_0F^{\nu}CH^r(\X)\subset F^{\nu}CH^r(\X)$ in the following way:
$$Z_0F^{\nu}CH^r(\X)= \sum_{\mathcal{Y},\G} Im\{ \G_{\ast}: F^{\nu}CH^{d_Y}(\mathcal{Y}) \la CH^r(\X)\},$$
where $\mathcal{Y}$ ranges over all projective smooth varieties with relative dimension $d_Y$ over $S$ and $\G$ runs over $CH^{r+l}(\mathcal{Y} \times \X)$.
Note that for $\nu=1$ we have $$Z_0F^1CH^r(\X)=CH^r(\X)_{alg}$$
where $CH^r(\X)_{alg}$ is the subgroup of cycle classes which are algebraically equivalent to zero.

\begin{defi}
We define the \textit{higher Griffiths group}
\begin{equation}\label{griff}Griff^{r,\nu}(\X)=F^{\nu}CH^r(\X)/F^{\nu+1}CH^r(\X)+Z_0F^{\nu}CH^r(\X).\end{equation}
\end{defi}
$Griff^{r,\nu}(\X)$ is a generalization of the Griffiths group in the context of the filtration on Chow groups. In fact,
for $\nu=1$, $Griff^{r,1}(\X)$ is the quotient of the Griffiths group $CH^r(\X)_{hom}/CH^r(\X)_{alg}$ by the image of $F^2CH^r(\X)$.

        \subsection{The infinitesimal invariant}

Let $\mathcal{C}$ be a family of smooth curves of genus $g$ over a smooth variety $S$ of dimension $n$, with a section $p$. Consider $\J$ the family of jacobian fibrations of relative dimension $g$ over $S$. Let $$ \C \la \J$$ be the canonical morphism of $\C$ into $\J$ defined by the section,
$P \la [P-p]$.
The image of this morphism is an algebraic cycle $\W$ of codimension $g-1$ in $\J$.
By the cycle class map, we take $\W \in CH^{g-1}(\J) \otimes \Q$. The Beauville's decomposition \cite{beau} allows us to decompose $CH^{g-1}(\J) \otimes  \Q$ in this way
$$\CH \otimes \Q= \bigoplus_{i=-1}^{g-1} CH^{g-1}_i(\J)_{\Q},$$
with
$ CH^{g-1}_i(\J)_{\Q}= \{ \alpha \in CH^{g-1}(\J) \otimes \Q \;| \;\;k^* \alpha= k^{2g-2-i}\alpha \;\;\forall k \in \Z \}$ where we denote with $k$ the multiplication on $\J$ and with $k^*$ the corresponding operation on the Chow group of $\J$.
So the cycle $\W$ is decomposed in the rational Chow group of $\J$ in the following way
$$ \W= \sum_{\nu } \W^{\nu} \in \CH \otimes \Q.$$
Moreover, Murre \cite{Murre} shows that
$\bigoplus_{i\geq \nu} CH^{g-1}_i(\J)_{\Q} \subset F^{\nu}\CH \otimes \Q.$\\

Consider the algebraic cycle $\W^{\nu} \in F^{\nu}\CH \otimes \Q$: for a general point $s \in S$, we define the \textit{higher infinitesimal invariant} as the linear map
$$\phi^{\nu}_s=\phi^{\nu}_s(\W^{\nu}):V^{g-1,\nu}_s \la \co$$
where $V^{g-1,\nu}_s$ is the cohomology of the Koszul complex
$$\bigwedge^{\nu+1}\tang \otimes H^0(\Om^{\nu+2}_{J_s}) \la \bigwedge^{\nu}\tang \otimes H^1(\Om^{\nu+1}_{J_s}) \la  \tang \otimes H^2( \Om^{\nu}_{J_s}).$$
Consider the natural map $\epsilon: T_{S,s} \tens H^0(\Om^{\nu+2}_{J_s}) \la H^1(\Om^{\nu+1}_{J_s})$ induced by the cup product with the Kodaira-Spencer class and
let $L_0V^{g-1,\nu}$ be the space given by $(\bigwedge^{\nu} T_{S,s} \tens Im (\epsilon)) \cap V^{g-1,\nu}_s$ (see \cite{Ikeda}).\\

\noindent
Ikeda (cfr. Prop. (3.7) in \cite{Ikeda}) shows a formula for calculating the infinitesimal invariant
\begin{equation}\label{ikeda formula} \phi^{\nu}_{s}(\W^{\nu})(\xi_1\wedge \cdots\wedge\xi_{\nu} \otimes w_1\wedge \cdots \wedge w_{\nu+1}\otimes \sigma)= < \alpha^{\nu+1}_{S,s}(w_1 \wedge \cdots \wedge w_{\nu+1}),\; \sigma >,\end{equation}
where $\xi_j \in T_{S,s}$, $w_i \in H^0(\Om^1_{J_s})$ and $\sigma \in H^1(\Ox_{J_s})$.
The computation of the infinitesimal invariant allows us to establish if the cycle is contained in $F^{\nu+1}CH^{g-1}(J_s)$ or in $Z_0F^{\nu}CH^{g-1}(J_s)$.
\begin{prop}[Ikeda] \label {prop ikeda}
Suppose $s \in S$ is generic with respect to $\W^{\nu} \in F^{\nu}\CH$.
\begin{enumerate}
\item If $\W^{\nu} \in F^{\nu+1}CH^{g-1}(J_s)$, then the infinitesimal invariant at $s$ $\phi^{\nu}_s(\W^{\nu}): V_s^{g-1,\nu} \la \co$ is zero.
\item If $\W^{\nu} \in Z_0F^{\nu}CH^{g-1}(J_s)$, then the infinitesimal invariant at $s$ $\phi^{\nu}_s(\W^{\nu}): L_0V_s^{g-1,\nu} \la \co$ is zero.
\end{enumerate}
\end{prop}
We refer to \cite{Ikeda} for the proof.\\

        \subsection{Main theorem}

In this section, we will prove \textbf{Theorem \ref{main thm intro 2}}. Consider a family of smooth curves $\mathcal{C}$ of genus $g\geq 10$ over a variety $S$ of dimension $2$ and the correspondent family of jacobian $\J$ over $S$. We have the following

\begin{thm}
For a family $\mathcal{C}$ of general smooth curves of genus $g \geq 16$ over a smooth variety $S$ of dimension $2$,
consider the algebraic cycle $\W^2 \in CH^{g-1}_2(\J)$ in the Beauville's decomposition of $\W$.
Then the cycle $\W_s^2$ is not trivial in $Griff^{g-1,2}(J_s)$ where $J_s$ is the fiber of the family $\J$ over the generic point $s \in S$. Moreover, there exists a family of curves of genus $g\geq 10$ such that the cycle $\W_s^2$ is a non trivial element in $Griff^{g-1,2}(J_s)$.

\end{thm}

\begin{proof}
Let $s\in S$ be a point. We compute the infinitesimal invariant at $s$ of the cycle $\W^2$.\\
Theorem \ref{main prop intro} assures if $C_s$ is a general smooth curve of genus $g \geq 16$, the adjunction map defined in (\ref{n=2}), $$\alpha^3_{S,s}: \bigwedge^3 U \la \mathcal{Q}_{\wc,\mathcal{V}}$$ is not zero, where $U$ is a $3$-dimensional subspace of $H^0(\Om^1_{\C}|_{C_s})$ contained in the kernel of the coboundary map of (\ref{omega}).
Moreover, again by Theorem \ref{main prop intro}, we know that there exists a family of curves of genus $g\geq 10$  for which the map $\alpha^3_{S,s}$ is not trivial.\\

Consider $w_i\in U$ for $i=1,2,3$. We remark that the condition $w_i \in U$ is equivalent to require that $$\xi_j \cdot w_i=0$$ for $i=1,2,3$ where $\xi_j \in H^1(T_{C_s})$ for $j=1,2$. Notice that $H^0(\Om^1_{C_s})\simeq H^0(\Om^1_{J_s})$; so we can take $w_i \in H^0(\Om^1_{J_s})$.
Take $w_4 \in H^0(\Om^1_{J_s})$ and set $\sigma= \xi_2 \cdot w_4 \in H^1(\mathcal{O}_{J_s})$. Obviously we have $<w_i, \bar \sigma >\,=0$ for $i=1,2,3$.\\

A direct computation shows that the element $\xi_1\wedge \xi_2 \tens w_1 \wedge w_2 \wedge w_3 \wedge \sigma$ is contained in $V^{g-1,\nu}_s$; the choice of $\sigma$ assures that $\xi_1\wedge \xi_2 \tens w_1 \wedge w_2 \wedge w_3 \wedge \sigma$ lies also in $L_0V^{g-1,2}_s$.\\

Therefore, we can evaluate $\phi^2_{s}(\xi_1\wedge \xi_2 \otimes w_1\wedge w_2 \wedge w_3\otimes \sigma).$
By formula (\ref{ikeda formula}) with $\nu=2$, we have $$\phi^2_{s}(\xi_1\wedge \xi_2 \otimes w_1\wedge w_2 \wedge w_3\otimes \sigma)
=<\alpha^3_{S,s}(w_1 \wedge w_2 \wedge w_3), \sigma>.$$  It is possible to verify directly this equality using the expression of $\alpha^3_{S,s}$ in coordinates. \\

\noindent
The non-triviality of $\alpha^3_{S,s}(w_1 \wedge w_2 \wedge w_3)$ allows us to conclude that also the infinitesimal invariant is not zero.
Proposition (\ref{prop ikeda}) concludes the proof. In fact, the cycle $\W_s^2$ is contained neither in $F^3CH^{g-1}(J_s)$ nor in $Z_0F^2CH^{g-1}(J_s)$. Hence $\W_s^2$ is a non trivial element in $Griff^{g-1,2}(J_s)$
where, by definition (\ref{griff}), we have $$Griff^{g-1,2}(J_s)= F^2CH^{g-1}(J_s) / F^{3}CH^{g-1}(J_s) + Z_0F^{2}CH^{g-1}(J_s).$$

\end{proof}

\vspace{0.4in}
\noindent
\textsc {Dipartimento di Matematica F. Casorati\\ Universit\`a di Pavia, via Ferrata 1, 27100 Pavia, Italy}\\
\textup {email: pirola@dimat.unipv.it}\\

\noindent
\textsc {Dipartimento di Matematica F. Brioschi\\ Politecnico di Milano, Piazza Leonardo da Vinci 32, 20133 Milano, Italy}\\
\textup {email: cecilia.rizzi@mate.polimi.it}
\end{document}